\documentclass[3p]{elsarticle}

\usepackage{amsmath}
\usepackage{amsthm}
\usepackage{amssymb}
\usepackage{bm}
\usepackage{accents}

\newtheorem{thm}{Theorem}
\newtheorem{lem}[thm]{Lemma}
\newdefinition{rmk}{Remark}
\newproof{pf}{Proof}

\newcommand{\const}{\mathop{\rm const}\nolimits}

\numberwithin{equation}{section}

\graphicspath{{figs/}}

\journal{arXiv}

\begin{document}

\begin{frontmatter}

\title{Numerical solution of the Cauchy problem for Volterra
integrodifferential equations with difference kernels\tnoteref{label1}}
\tnotetext[label1]{The publication was financially supported by the research grant 20-01-00207 of Russian Foundation for Basic Research}

\author{P.N. Vabishchevich\corref{cor1}\fnref{lab1,lab2}}
\ead{vabishchevich@gmail.com}
\cortext[cor1]{Correspondibg author.}

\address[lab1]{Nuclear Safety Institute, Russian Academy of Sciences,
              52, B. Tulskaya, 115191 Moscow, Russia}

\address[lab2]{North-Caucasus Center for Mathematical Research, North-Caucasus Federal University, 
                   1, Pushkin Street, 355017 Stavropol, Russia}

\begin{abstract}
We consider the problems of the numerical solution of the Cauchy problem for an evolutionary equation with memory when the kernel of the integral term is a difference one.
The computational implementation is associated with the need to work with an approximate solution for all previous points in time.
In this paper, the considered nonlocal problem is transformed into a local one; a loosely coupled equation system with additional ordinary differential equations is solved.
This approach is based on the approximation of the difference kernel by the sum of exponentials.
Estimates for the stability of the solution concerning the initial data and the right-hand side for the corresponding Cauchy problem are obtained.
Two-level schemes with weights with convenient computational implementation are constructed and investigated.
The theoretical consideration is supplemented by the results of the numerical solution of the integrodifferential equation when the kernel is the stretching exponential function.
\end{abstract}


\end{frontmatter}

\section{Introduction}\label{sec:1}

Classical applied mathematical models are based on systems of partial differential equations.
Nonstationary processes are described by parabolic and hyperbolic equations \cite{LionsBook,evans2010partial}.
More general models, which partially inherit both the properties of parabolic and hyperbolic equations,
are associated with evolutionary integrodifferential equations \cite{GripenbergBook1990,pruss2013evolutionary}.
The simplest nonlocal in time problems are generated with fractional time derivatives \cite{baleanu2012fractional,uchaikin}, approximate methods for which have been actively discussed over the last decade.

The numerical study of multidimensional boundary value problems for equations with memory is carried out with standard finite element or finite volume approximations in the space \cite{KnabnerAngermann2003,QuarteroniValli1994}.
The focus is on (see, for example, \cite{ChenBook1998}) the construction of time approximations.
The most straightforward approach \cite{mclean1993numerical,mclean1996discretization} is associated with using certain quadratures for the integral term and the usual two-level approximations of the time derivative:  backward Euler and Crank-Nicolson schemes.

The evolutionary equations for Volterra integrodifferential equations are not very computationally convenient.
The solution depends on the complete history, and therefore, when transformed to a new level in time, we need to work with the solution on all previous levels in time.
Basic approaches to reduce computational work for fractional time derivative problems are discussed in \cite{diethelm2020good}.
The most excellent possibilities are provided by methods based on kernel approximation, which leads us to more simpler problems.

We are interested in kernel approximations that can significantly reduce the computational work for an approximate solution of the integrodifferential problem. In particular, we want to reduce the memory requirement.
For the Volterra integral equations, a transformation approach to the set of Cauchy problems for a system of ODEs has long been known.
One of the first works in this direction of research is the article \cite{bownds1977note}.
The book \cite{linz1985analytical} separately discusses two essential cases for Volterra integral equations with a difference kernel.
The first one is related to the case when the kernel is a polynomial.
In the second, more interesting case, the kernel is the sum of exponents.
It is this approach that is used in this paper.

We formulate the Cauchy problem for the Volterra integrodifferential equation with a positive definite self-adjoint operator in a real finite-dimensional Hilbert space.
The sum of exponentials approximates the difference kernel with fixed sigh coefficients.
The nonlocal problem is transformed to the Cauchy problem for a local system of equations, including the standard evolution equation for the problem operator and an additional set of ordinary differential equations.
A priori estimates are given to solve this problem, consistent with the a priori estimate for the original evolutionary equation with a positive-type memory term.
The main result of the work consists of constructing two-level difference schemes and formulating conditions for the stability of an approximate solution. 

The paper is organized as follows.
Section \ref{sec:2} describes an integrodifferential equation with a positive definite operator in a finite-dimensional Hilbert space.
An estimate of the stability of the solution for problems with a positive definite difference kernel is given.
Kernel approximation by the sum of exponentials is discussed in Section \ref{sec:3}.
In Section \ref{sec:4}, we formulate a system of local equations for an approximate solution to a nonlocal problem.
Estimates of the stability of the solution based on the initial data and the right-hand side are obtained.
We build unconditionally stable two-level schemes in Section \ref{sec:5}.
Section \ref{sec:6} is devoted to more general evolutionary memory problems.
Numerical experiments for a two-dimensional model problem for integrodifferential equation are given in Section \ref{sec:7}.
In our calculations, the kernel is the stretching exponential function.
The results are summarized in Section \ref{sec:8}.

\section{Problem formulation}\label{sec:2}

In a bounded polygonal domain $\Omega \subset \mathbf{R}^d, \ d = 2,3$ with the Lipschitz continuous boundary $\partial\Omega$
we introduce the elliptic operator as
\begin{equation}\label{2.1}
  \mathcal{A}  v = - {\rm div}  ( a({\bm x}) {\rm grad} \, v) + c({\bm x}) v
\end{equation} 
with coefficients $0 < a_1 \leq a({\bm x}) \leq a_2$, $c({\bm x}) \geq 0$.
The operator $\mathcal{A}$ is defined on the set of functions $v({\bm x})$ that satisfy
on the boundary $\partial\Omega$ the following conditions:
\begin{equation}\label{2.2}
  v ({\bm x}) = 0,
  \quad {\bm x} \in \partial \Omega .
\end{equation} 
In the Hilbert space $\mathcal{H} = L_2(\Omega)$, the operator $\mathcal{A}: \mathcal{H} \rightarrow \mathcal{H}$ is self-adjoint and positive definite:
\begin{equation}\label{2.3}
  \mathcal{A}  = \mathcal{A}^* \geq \nu_{\mathcal{A}}   \mathcal{I} ,
  \quad \nu_{\mathcal{A}}  > 0 ,    
\end{equation} 
where $\mathcal{I}$ is the identity operator in $\mathcal{H}$.

The Cauchy problem for an evolutionary equation with memory is considered.
The function $w(\bm x,t)$ ($w(t) = w(\cdot, t)$) satisfies an integrodifferential equation with a difference kernel
\begin{equation}\label{2.4}
 \frac{d w}{d t} + \int_{0}^{t} k(t-s)  \mathcal{A} w(s) d s = f(t),
 \quad t > 0 ,  
\end{equation} 
and the initial condition
\begin{equation}\label{2.5}
 w(0) = w_0 .
\end{equation} 

The kernel $k(t)$  (see for example \cite{mclean1993numerical,mclean1996discretization}) is assumed to be real-valued and
positive-definite, i.e., for each $T > 0$ the kernel $k(t)$ belongs to $L_1(0, T)$ and satisfies
\begin{equation}\label{2.6}
 \int_{0}^{T} \psi(t) \int_{0}^{t} k(t-s) \psi(s) d s \, d t \geq 0,
 \quad \psi \in C[0,T] .
\end{equation} 
These kernels are often called convolution kernels of positive type \cite{gripenberg1990volterra}.
The sufficient condition for the positive-definite kernel $k(t)$ \cite{halanay1965asymptotic} is 
\begin{equation}\label{2.7}
 k(t) \geq 0,
 \quad \frac{d k}{d t}  (t) \leq 0,
 \quad \frac{d^2 k}{d t^2}  (t) \geq 0,
 \quad t > 0 . 
\end{equation} 

With finite element or finite volume approximation in space \cite{KnabnerAngermann2003, QuarteroniValli1994}, the operator
$\mathcal {A}$ is associated with an operator $A$ in a finite-dimensional Hilbert space $H$.
In $H$, we use the usual notation $(\cdot, \cdot)$ and $\|\cdot\|$ to the scalar product and norm.
For $A: H \rightarrow H$ the property of self-adjointness and positive definiteness (\ref{2.3}) is preserved:
\begin{equation}\label{2.8}
  A = A^* \geq \nu_A  I,
  \quad \nu_A   > 0 ,    
\end{equation} 
where $I$ is the identity operator in $H$.

The Cauchy problem (\ref{2.4}), (\ref{2.5}) corresponds to the problem
\begin{equation}\label{2.9}
 \frac{d u}{d t} + \int_{0}^{t} k(t-s)  A u(s) d s = \varphi (t),
 \quad t > 0 ,  
\end{equation} 
\begin{equation}\label{2.10}
 u(0) = u_0 .
\end{equation} 
The convergence of the approximate solution $ u (t) $ to $ w (t) $ with finite element approximation in space is proved in \cite{mclean1993numerical}.
It is based on the corresponding a priori estimates for solving the problem (\ref{2.9}), (\ref{2.10}).
In our study, the guideline is the following statement.

\begin{thm}\label{t-1}
Let $A$  be a positive definite self-adjoint operator in a Hilbert space and let $k(t)$ be a positive definite kernel. 
Then for the solution of (\ref{2.9}), (\ref{2.10}) we have the estimate of stability of the initial data and the right-hand side
\begin{equation}\label{2.11}
 \|u(t)\| \leq \|u_0\| + \int_{0}^{t} \|\varphi(s)\| d s,
 \quad t > 0 . 
\end{equation} 
\end{thm}

\begin{pf}
We multiply the equation (\ref{2.9}) scalarly in $H$ by $u(t)$ and obtain
\[
 \frac{1}{2} \frac{d}{d t} \|u(t)\|^2 + \int_{0}^{t} k(t-s)  \big(A u(s), u(t) \big) d s = (\varphi (t), u(t)) .
\] 
By integration over $(0, T)$ this yields 
\[
 \frac{1}{2} \big(\|u(T)\|^2 - \|u(0)\|^2 \big) +  \int_{0}^{T} \int_{0}^{t} k(t-s)  \big(A u(s), u(t) \big) d s \ d t 
 \leq  \int_{0}^{T}  \|\varphi (t)\| \|u(t))\| d t.
\]  
Taking into account properties (\ref{2.6}), (\ref{2.8}) and the initial condition (\ref{2.10}), we have the inequality
\begin{equation}\label{2.12}
 \|u(T)\|^2 \leq \|u_0\|^2 + 2 \int_{0}^{T}  \|\varphi (t)\| \|u(t))\| d t.
\end{equation} 
We define
\[
 R(T) = \int_{0}^{T}  \|\varphi (t)\| \|u(t))\| d t ,
\] 
therefore
\[
 \frac{d R}{d T} = \|\varphi (T)\| \|u(T))\| .
\] 
Since (\ref{2.12}), we have
\[
 \frac{d R}{d T} \leq \|\varphi (T)\| \big(\|u_0\|^2 + 2 R(T)\big)^{1/2} .
\]
Whereas
\[
 \frac{d }{d T} \big(\|u_0\|^2 + 2 R(T)\big)^{1/2} =
 \big(\|u_0\|^2 + 2 R(T)\big)^{-1/2} \frac{d R}{d T} ,
\] 
we get
\[
 \frac{d }{d T} \big(\|u_0\|^2 + 2 R(T)\big)^{1/2} \leq \|\varphi (T)\| . 
\] 
Integration of this inequality with respect to $T$ from $0$ to $t$ with the inequality (\ref{2.12}) and $R(0) = 0$ gives the desired estimate (\ref{2.11}).
\end{pf}

\begin{rmk}
We have given a direct proof that the linear inequality (\ref{2.11}) follows from the nonlinear inequality (\ref{2.12}) for $\|u(T)\|$.
We could also come to this result if we used the Bihari lemma \cite{bihari1956generalization,beckenbach2012inequalities}
or a more particular result --- the Vaigant inequality \cite{kaliev2006boundary,qin2017analytic}. 
\end{rmk}

\begin{rmk}
For a self-adjoint and positive operator $D$, the Hilbert space $H_D$ is defined
with scalar product and norm $(u, v)_D = (D u, v), \ \| u \|_D = (u, v)_D ^ {1/2} $.
For operators $D$ that commute with $A$, similarly (\ref{2.11}), a more general estimate 
\begin{equation}\label{2.13}
 \|u(t)\|_D \leq \|u_0\|_D + \int_{0}^{t} \|\varphi(s)\|_D d s,
 \quad t > 0 , 
\end{equation} 
is established.
\end{rmk}

The standard approach to approximation (\ref{2.9}) in time is based on choosing a specific difference derivative for the first term of the equation and approximating the integral expression using the corresponding quadrature formula.
On this path, we are building, for example, analogs of the classic the backward Euler and Crank-Nicolson schemes \cite{mclean1993numerical,mclean1996discretization}.
The main disadvantage of such methods is the need to remember the solution at previous points in time.
Various approaches are used for the treatment of the memory term in nonlocal problems.
For fractional derivative problems, see \cite{diethelm2020good} for a general discussion of this problem.
Of most significant interest have approaches with constructing an approximate solution to a nonlocal problem with memory using solutions of local problems.

\section{Approximation of the kernel by the sum of exponentials}\label{sec:3}

In the approximation of nonlinear functions, rational approximations are most widely used in \cite{braess1986nonlinear}.
A separate class of methods for approximating functions is associated with approximations of the sum of exponents.
When approximating the kernel of integrodifferential equations with difference kernel, the study is based on applying the Laplace transform.
In this case, the rational approximation of the Laplace transform of the kernel is directly related to the approximations by the sum of exponentials for the kernel.
 
The kernel $k(t)$ is approximated by the function $\accentset {\sim}{k}(t)$, which has the form (Prony series)
\begin{equation}\label{3.1}
 \accentset{\sim }{k}(t) = \sum_{i=1}^{m} a_i \exp(-b_i t),
 \quad t \geq 0 . 
\end{equation} 
For the coefficients $a_i, b_i, \ i = 1,2,\ldots, m,$ the assumptions are fulfilled:
\begin{equation}\label{3.2}
 a_i > ,
 \quad b_i \geq  0,
 \quad   i = 1,2,\ldots, m .
\end{equation} 
Taking into account the conditions (\ref{2.7}) under the constraints (\ref{3.2}), the kernel $\accentset {\sim}{k}(t) $ is a positive definite.

An approximate solution to the problem (\ref{2.9}), (\ref{2.10}) is denoted by $v(t)$; it is defined as a solution to the Cauchy problem
\begin{equation}\label{3.3}
 \frac{d v}{d t} + \int_{0}^{t} \accentset{\sim }{k}(t-s)  A v(s) d s = \varphi (t),
 \quad t > 0 ,  
\end{equation} 
\begin{equation}\label{3.4}
 v(0) = u_0 .
\end{equation} 
An analogue of (\ref{2.11}) is the estimate
\begin{equation}\label{3.5}
 \|v(t)\| \leq \|u_0\| + \int_{0}^{t} \|\varphi(s)\| d s,
 \quad t > 0 . 
\end{equation} 
The accuracy of the approximate solution is determined by the accuracy of the kernel approximation.

\begin{thm}\label{t-2}
Let the kernel approximation error satisfy the estimate
\begin{equation}\label{3.6}
 |\accentset{\sim }{k}(t) - k(t)| \leq \varepsilon ,
 \quad t \geq 0 .  
\end{equation} 
Then, for the error of the approximate solution $v - u$, the estimate
\begin{equation}\label{3.7}
 \|v(t)-u(t)\| \leq \varepsilon \int_{0}^{t} s \Big ( \|Au_0\| + \int_{0}^{s} \|A \varphi (\theta) \| d \theta \Big ) d s,
 \quad t > 0 , 
\end{equation} 
is holds.
\end{thm}

\begin{pf}
From Cauchy problems (\ref{2.9}), (\ref{2.10}) and (\ref{3.3}), (\ref{3.4}) for $w = v-u$ we obtain
\begin{equation}\label{3.8}
 \frac{d w}{d t} + \int_{0}^{t} k(t-s)  A w(s) d s = 
 \int_{0}^{t} \big(k(t-s) - \accentset{\sim }{k}(t-s)\big)   A v(s) d s,
 \quad t > 0 ,  
\end{equation} 
\begin{equation}\label{3.9}
 w(0) = 0 .
\end{equation} 
For a solution of the problem (\ref{3.8}), (\ref{3.9}), the following estimate holds:
\begin{equation}\label{3.10}
 \|w(t)\| \leq \int_{0}^{t} \|\psi(s)\| d s,
 \quad t > 0 , 
\end{equation} 
wherein
\[
 \psi(t) = \int_{0}^{t} \big(k(t-s) - \accentset{\sim }{k}(t-s)\big)   A v(s) d s .
\] 
Setting  $D = A^2$ (see (\ref{2.13})) and taking into account the estimate (\ref{3.6}) we obtain
\[
 \|\psi(s)\| \leq \varepsilon s \Big ( \|Au_0\| + \int_{0}^{s} \|A \varphi (\theta) \| d \theta \Big ) .
\] 
Under these conditions, the error estimate (\ref{3.7}) follows from the inequality (\ref{3.10}).
\end{pf}

\section{System of local equations}\label{sec:4} 

The use of the kernel $k(t)$ approximations by the sum of exponentials allows one to go from a nonlocal equation (\ref{2.9}) to a system of local equations \cite{linz1985analytical}.
This approach for fractional derivative ODEs is known as kernel compression schemes \cite{diethelm2020good,baffet2017kernel,baffet2017high}.

Consider the problem  (\ref{2.9}), (\ref{2.10}) in the special case when
\[
 k(t) = a_1\exp(-b_1 t),
 \quad a_1 > 0, \quad b_1 > 0 .   
\] 
Differentiation of equation
\[
 \frac{d u}{d t} + a_1\int_{0}^{t} \exp(-b_1 (t-s))  A u(s) d s = \varphi (t),
\] 
gives us
\begin{equation}\label{4.1}
  \frac{d^2 u}{d t^2} + b_1 \frac{d u}{d t} + a_1 A u = \frac{d \varphi}{d t} + b_1 \varphi (t) .  
\end{equation} 
For a local second-order evolution equation (\ref{4.1}), the Cauchy problem is posed and the initial conditions are of the form
\begin{equation}\label{4.2}
 u(0) = u_0,
 \quad  \frac{d u}{d t}(0) = \varphi (0) .
\end{equation} 
We would like to have a similar transformation from a nonlocal problem to a local one for a more general case (\ref{3.1})--(\ref{3.4}).  

The problem of (\ref{4.1}), (\ref{4.2}) is not very useful for the study.
The standard a priori estimate for a second-order evolution equation is obtained by multiplying the equation by $d u / dt$.
Applied to the equation (\ref{4.1}), this gives the inequality
\[
 \frac{d}{d t} \Big( \Big \|\frac{d u}{d t} \Big\|^2 + a_1 \|u\|_A^2 \Big) \leq \frac{b_1}{2}  \|\eta (t)\|^2,
 \quad \eta (t) = \frac{d \varphi}{d t} + b_1 \varphi (t) .  
\] 
The thus obtained estimate
\[
 \Big \|\frac{d u}{d t}(t) \Big\|^2 + a_1 \|u(t)\|_A^2 \leq \|\varphi (0)\|^2 + a_1 \|u_0\|_A^2 + \frac{b_1}{2} \int_{0}^{t}  \|\eta (s)\|^2 d s
\] 
difficult to directly relate to the estimate (\ref{3.5}) for the problem (\ref{4.1}), (\ref{4.2}).

Instead of a second-order evolution equation (\ref{4.1}), we will consider a system of two first-order evolution equations.
We put
\[
 u_1(t) = \int_{0}^{t} \exp(-b_1 (t-s)) u(s) d s ,
\] 
which gives the system of equations
\begin{equation}\label{4.3}
 \frac{d u}{d t} + a_1 A u_1 = \varphi (t) , 
\end{equation} 
\begin{equation}\label{4.4}
 \frac{d u_1}{d t} + b_1 u_1 - u = 0 .
\end{equation} 
The initial conditions are
\begin{equation}\label{4.5}
 u(0) = u_0,
 \quad u_1(0) = 0 . 
\end{equation} 
Multiplying equation (\ref{4.3}) to $u(t)$ and equation (\ref{4.4}) to $a_1 A u_1 (t)$, we obtain
\[
 \frac{1}{2}  \frac{d}{d t} \|u\|^2 + a_1 (A u,u_1) = (\varphi, u),
\] 
\[
 \frac{a_1}{2} \frac{d}{d t} \|u_1\|_{A}^2 + a_1 b_1  \|u_1\|_{A}^2 - a_1 (A u,u_1) = 0 .
\] 
Adding these relations and taking into account the initial conditions (\ref{4.5}) and the positiveness of $a_1, b_1$, we obtain
\begin{equation}\label{4.6}
 \Big (\|u(t)\|^2 + a_1 \|u_1(t)\|_{A}^2 \Big)^{1/2} \leq \|u_0\| +  \int_{0}^{t} \|\varphi(s)\| d s .
\end{equation} 
The estimate (\ref{3.5}) follows from (\ref{4.6}).

The Cauchy problem (\ref{3.1})--(\ref{3.4}) is treated similarly.
Let us introduce the functions
\[
 v_i(t) = \int_{0}^{t} \exp(-b_i (t-s)) u(s) d s ,
 \quad i = 1,2, \ldots, m . 
\] 
Given (\ref{3.1}) the equation (\ref{3.3}) can be written as
\begin{equation}\label{4.7}
 \frac{d v}{d t} + \sum_{i=1}^{m} a_i A v_i(t) = \varphi (t) .
\end{equation} 
Similarly the equation (\ref{4.4}), we have
\begin{equation}\label{4.8}
 \frac{d v_i}{d t} + b_i v_i - v = 0 ,
 \quad i = 1,2, \ldots, m . 
\end{equation} 
The system of equations (\ref{4.7}), (\ref{4.8}) is supplemented with the initial conditions
\begin{equation}\label{4.9}
 v(0) = u_0,
 \quad v_i(0) = 0,
 \quad i = 1,2, \ldots, m .  
\end{equation} 

\begin{thm}\label{t-3}
Let $A$  be a positive definite self-adjoint operator in a Hilbert space $H$. 
Then for the solution to the problem (\ref{3.2}), (\ref{4.7})--(\ref{4.9}), the stability estimate of the solution with respect to the initial data and the right-hand side 
\begin{equation}\label{4.10}
 \Big (\|v(t)\|^2 + \sum_{i=1}^{m} a_i \|v_i(t)\|_{A}^2 \Big)^{1/2} \leq \|u_0\| +  \int_{0}^{t} \|\varphi(s)\| d s ,
 \quad t > 0 , 
\end{equation} 
is valid.
\end{thm}

\begin{pf}
As in the proof the estimate (\ref{4.6}), we multiply the equation (\ref{4.7}) by $v(t)$ and the separate equation (\ref{4.8}) for $i = 1,2, \ldots, m$ to $a_i A v_i (t)$. This gives equalities
\[
 \frac{1}{2}  \frac{d}{d t} \|v\|^2 + \sum_{i=1}^{m} a_i (A v,v_i) = (\varphi, v),
\] 
\[
 \frac{a_i}{2} \frac{d}{d t} \|v_i\|_{A}^2 + a_i b_i  \|v_i\|_{A}^2 - a_i (A v,v_i) = 0 ,
 \quad i = 1,2, \ldots, m .  
\] 
Adding them taking into account the initial conditions (\ref{4.9}) and constraints (\ref{3.2}) on the coefficients $a_i, b_i, \  i = 1,2, \ldots, m$, we obtain
\[
 \frac{1}{2}  \frac{d}{d t} \Big ( \|v\|^2 + \sum_{i=1}^{m} a_i \|v_i\|_{A}^2 \Big ) \leq \|\varphi \| \| v\| .
\] 
It follows proves estimate (\ref{4.10}).
\end{pf}

It is convenient to write the system (\ref{4.7}), (\ref{4.8}) in the form of one first-order equation for vector quantities.
Define a vector 
$\bm v = \{v, v_1, \ldots, v_m \} $ and $\bm \varphi  = \{\varphi, 0, \ldots, 0 \} $
and from (\ref{4.7}), (\ref{4.8}), we get to the Cauchy problem
\begin{equation}\label{4.11}
 \bm B \frac{d \bm v}{d t} + \bm A \bm v = \bm \varphi  ,
\end{equation} 
\begin{equation}\label{4.12}
 \bm v(0) = \bm v_0 ,
\end{equation} 
where $\bm v_0 = \{u_0, 0, \ldots, 0 \}$.
For the operator matrices $\bm B$ and $\bm A$,  we have the representation
\[
 \bm B = \mathrm{diag} \, \big (I, a_1 A, \ldots, a_m A \big),
 \quad \bm A = \left (\begin{array}{cccc}
  0  &  a_1 A & \cdots &  a_m A \\
  - a_1 A  &  a_1 b_1 A & \cdots &  0 \\
  \cdots  & \cdots & \cdots &  0 \\
  - a_m A &  0 & \cdots &  a_m b_m A  \\
\end{array}
 \right ) .
\]

Problem (\ref{4.11}), (\ref{4.12}) can be 
considered on the direct sum of spaces $\bm H = H \oplus \ldots  \oplus H$,
when for $\bm v, \bm w \in \bm H$,  the scalar product and norm are determined by the expressions
\[
 (\bm v, \bm w) = (v, w) + \sum_{1=1}^{m} (v_i, w_i),
 \quad \|\bm v\| =  (\bm v, \bm v)^{1/2} .
\]
Taking into account properties (\ref{2.8}) and (\ref{3.2}), we have
\begin{equation}\label{4.13}
 \bm B = \bm B^* > 0,
 \quad \bm A \geq 0 . 
\end{equation} 

To prove the estimate (\ref{4.10}), we multiply the equation (\ref{4.11}) scalarly in $\bm H$ by $\bm v$.
Given (\ref{4.13}) it follows
\[
 \|\bm v\|_{\bm B} \frac{d}{d t}  \|\bm v\|_{\bm B} \leq (\bm \varphi , \bm v) .  
\]
Since
\[
 (\bm \varphi , \bm v) \leq \|\bm \varphi \|_{\bm B^{-1}} \|\bm v\|_{\bm B} 
\] 
we have the standard a priori estimate for the Cauchy problem (\ref{4.11}), (\ref{4.12}) 
\begin{equation}\label{4.14}
 \|\bm v(t)\|_{\bm B} \leq \|\bm v_0\|_{\bm B} + \int_{0}^{t} \|\bm \varphi(s)\|_{\bm B^{-1}} d s .
\end{equation} 
In our case
\[
 \|\bm v(t)\|_{\bm B} = \Big (\|v(t)\|^2 + \sum_{i=1}^{m} a_i \|v_i(t)\|_{A}^2 \Big)^{1/2} , 
 \quad \|\bm v_0\|_{\bm B} = \|u_0\| ,
 \quad \|\bm \varphi(t)\|_{\bm B^{-1}} = \|\varphi(t)\|,
\] 
therefore the estimate (\ref{4.14}) gives (\ref{4.10}).

\section{Time approximation}\label{sec:5} 

In the approximate solution of the Cauchy problem (\ref{4.11}), (\ref{4.12}),  implicit time approximations are often used.
In this case, we have unconditionally stable schemes. 
We will use, for simplicity, a uniform grid in time with step $\tau$ and let  $y^n=y(t^n), \ t^n = n\tau$, $n =0, 1, \ldots$. 
We will use a two-level scheme with the weight $\sigma = \const \in (0,1]$, when 
\begin{equation}\label{5.1}
 \bm B \frac{\bm y^{n+1} - \bm y^{n}}{\tau } + \bm A \bm y^{n+\sigma} = \bm \varphi^{n+\sigma},
 \quad n = 0,1,\ldots,
\end{equation} 
\begin{equation}\label{5.2}
 \bm y^0 = \bm v_0 .
\end{equation} 
We used the notation
\[
 \bm y^{n+\sigma} = \sigma \bm y^{n+1} + (1-\sigma ) \bm y^{n} .
\] 

The difference scheme (\ref{5.1}), (\ref{5.2}) approximates the Cauchy problem (\ref{4.11}), (\ref{4.12}) with sufficient smoothness of the solution $\bm v(t)$ with the first order in $\tau$ for $\sigma \neq 0.5$ and with the second order for $\sigma = 0.5$ (Crank-Nicolson scheme).
To study the stability of two-level schemes, we can use the results of the theory of stability (correctness) of operator-difference schemes \cite{SamarskiiTheory,SamarskiiMatusVabischevich2002}.

\begin{thm}\label{t-4}
The two-level scheme (\ref{4.13}), (\ref{5.1}), (\ref{5.2}) 
is unconditionally stable for $\sigma \geq 0.5$.
Under these constraints, for an approximate solution to the problem (\ref{4.11}), (\ref{4.12}), the a priori estimate
\begin{equation}\label{5.3}
 \|\bm y^{n+1}\|_{\bm B} \leq \|\bm v_0\|_{\bm B}  + \sum_{k=0}^{n}\tau  \|\bm \varphi^{k+\sigma}\|_{\bm B^{-1}} ,
 \quad n = 0,1,\ldots ,
\end{equation} 
holds.
\end{thm}

\begin{pf}
Multiplying the equation (\ref{5.1}) to $\tau \bm y^{n+\sigma}$, we get the inequality
\begin{equation}\label{5.4}
 (\bm B (\bm y^{n+1} - \bm y^{n}), \bm y^{n+\sigma}) \leq \tau \|\bm \varphi^{n+\sigma}\|_{\bm B^{-1}}  \|\bm y^{n+\sigma }\|_{\bm B} .
\end{equation}  
To evaluate the left-hand side of the inequality (\ref{5.4}), the following result is used \cite{vabishchevich2013flux}.

\begin{lem}\label{l-1}
Let be
\[
  w = \sigma u + (1-\sigma) v  
\] 
for $u, v$ from some Hilbert space $H_D$ ($D=D^* >0$).
Then
\begin{equation}\label{5.5}
  (D(u-v),w) \geq (\|u\|_D - \|v\|_D) \|w\|_D .
\end{equation}
if constant $\sigma \geq 0.5$. 
\end{lem}

Taking into account the inequality (\ref{5.5}), we get
\[
 \big (\bm B (\bm y^{n+1} - \bm y^{n}), y^{n+\sigma} \big ) \geq 
 \big (\|\bm y^{n+1}\|_{\bm B} - \|\bm y^{n}\|_{\bm B} \big ) \|\bm y^{n+\sigma}\|_{\bm B} .
\] 
With this in mind, from the inequality (\ref{5.4}) it follows
\[
 \|\bm y^{n+1}\|_{\bm B} \leq  \|\bm y^{n}\|_{\bm B} + \tau \|\bm \varphi^{n+\sigma}\|_{\bm B^{-1}} .
\] 
Thus, the a priori estimate (\ref{5.3}) holds, which acts as a grid analogue of the estimate (\ref{4.14}).
\end{pf}

Let's write a scheme with weights (\ref{5.1}), (\ref{5.2}) for individual components.
For an approximate solution of the problem (\ref{4.7})--(\ref{4.9}), the difference scheme is used
\begin{equation}\label{5.6}
 \frac{y^{n+1} - y^{n}}{\tau } + \sum_{i=1}^{m} a_i A y_i^{n+\sigma} = \varphi^{n+\sigma },
\end{equation} 
\begin{equation}\label{5.7}
 \frac{y_i^{n+1} - y_i^{n}}{\tau} + b_i y_i^{n+\sigma } - y^{n+\sigma} = 0 ,
 \quad i = 1,2, \ldots, m ,
 \quad n = 0,1, \ldots ,  
\end{equation} 
\begin{equation}\label{5.8}
 y^{0} = u_0,
 \quad y_i^{0} = 0,
 \quad i = 1,2, \ldots, m .  
\end{equation} 
From the estimate (\ref{5.3}) follows an a priori estimate
\begin{equation}\label{5.9}
 \Big (\|y^{n+1}\|^2 + \sum_{i=1}^{m} a_i \|y_i^{n+1}\|_{A}^2 \Big)^{1/2} \leq \|u_0\| + \sum_{k=0}^{n}\tau \|\varphi^{n+\sigma }\| ,
 \quad n = 0,1, \ldots ,
\end{equation} 
for the solution of the problem (\ref{5.6})--(\ref{5.8}).
The estimate (\ref{5.9}) is a difference analogue of the estimate (\ref{4.10}) for solving the differential problem (\ref{4.7})--(\ref{4.9}).

The problem of computational implementation deserves special attention when solving nonlocal problems.
In the case the scheme (\ref{5.6})--(\ref{5.7}) from the equation (\ref{5.7}), we have
\begin{equation}\label{5.10}
 y_i^{n+1} = \frac{\sigma \tau }{1 + \sigma b_i \tau } y^{n+1} + \chi_i^{n} ,
 \quad \chi_i^{n} =  \frac{1 }{1 + \sigma b_i \tau } \Big ( (1-\sigma) \tau y^{n} + \big (1 - (1-\sigma) b_i \tau ) \big ) y_i^{n} \Big ) ,
 \quad i = 1,2,\ldots, m .
\end{equation}  
Substitution in the equation (\ref{5.6}) gives the equation for finding $y^{n + 1}$
\begin{equation}\label{5.11}
 (1 + \mu \sigma \tau  A) y^{n+1} = \chi^{n} .
\end{equation} 
For the coefficient $\mu$ and the right-hand side, we have
\[
 \mu = \sum_{i=1}^{m} \frac{\sigma a_i \tau }{1 + \sigma b_i \tau } ,
 \quad \chi^{n} = y^{n} + \tau \varphi^{n+\sigma } - \tau \sum_{i=1}^{m} a_i A \big ( (1-\sigma) y_i^{n} + \sigma \chi_i^{n} \big ) .
\] 
Thus, the transition to a new $n+1$ level in time is provided by solving the standard problem (\ref{5.11}) for $y^{n+1}$ and calculating the auxiliary values $y_i^{n+1}, \ i = 1,2,\ldots, m,$ according to (\ref{5.11}).
The computational complexity of an approximate solution to the considered nonlocal problem (\ref{2.1}), (\ref{2.4}), (\ref{2.5}) is not much greater than when solving a local evolutionary problem with the operator $A$. 
It is necessary to additionally operate with solutions $m$ of simple auxiliary local evolutionary problems when explicitly calculating their solutions on a new level in time.

\section{More general nonlocal problems}\label{sec:6} 

The considered transformation of a nonlocal problem to a local one can be applied for more general than (\ref{2.8})--(\ref{2.10}) problems.
We will assume that instead of the equation (\ref{2.9}), we use the operator integrodifferential equation
\begin{equation}\label{6.1}
 B \frac{d u}{d t} + \int_{0}^{t} k(t-s)  A u(s) d s + C u = \varphi (t),
 \quad t > 0 ,  
\end{equation} 
in which, in addition to the property (\ref{2.8}), 
\begin{equation}\label{6.2}
 B = B^* \geq \delta_B I , 
 \quad \delta_B > 0 ,
 \quad C \geq 0 .  
\end{equation} 
When the kernel is approximated according to (\ref{3.1}), (\ref{3.2}), we come from the equation (\ref{6.1}) to the equation
\begin{equation}\label{6.3}
 B \frac{d v}{d t} + \int_{0}^{t} \accentset{\sim }{k}(t-s)  A v(s) d s + C v = \varphi (t),
 \quad t > 0 .  
\end{equation} 

When using auxiliary functions $v_i(t), \ i = 1,2, \ldots, m,$ the equation (\ref{6.3}) is written as
\begin{equation}\label{6.4}
 B \frac{d v}{d t} + \sum_{i=1}^{m} a_i A v_i(t) + C v = \varphi (t) .
\end{equation} 
The system of equations (\ref{4.7}), (\ref{6.4}) has the form (\ref{4.11}), where
\[
 \bm B = \mathrm{diag} \, \big (B, a_1 I, \ldots, a_m I \big),
 \quad \bm A = \left (\begin{array}{cccc}
  C  &  a_1 A & \cdots &  a_m A \\
  - a_1 A  &  a_1 b_1 A & \cdots &  0 \\
  \cdots  & \cdots & \cdots &  0 \\
  - a_m A &  0 & \cdots &  a_m b_m A  \\
\end{array}
 \right ) .
\]
Taking into account properties (\ref{2.8}) and (\ref{6.2}), conditions (\ref{4.13}) are valid.
In the estimate (\ref{4.14}), we now have
\[
 \|\bm v(t)\|_{\bm B} = \Big (\|v(t)\|_B^2 + \sum_{i=1}^{m} a_i \|v_i(t)\|_{A}^2 \Big)^{1/2} , 
 \quad \|\bm v_0\|_{\bm B} = \|u_0\|_B ,
 \quad \|\bm \varphi(t)\|_{\bm B^{-1}} = \|\varphi(t)\|_{B^{-1}} .
\] 
For a solution of the Cauchy problem (\ref{4.8}), (\ref{4.9}), (\ref{6.4}), this gives the estimate
\begin{equation}\label{6.5}
 \Big (\|v(t)\|^2 + \sum_{i=1}^{m} a_i \|v_i(t)\|_{A}^2 \Big)^{1/2} \leq \|u_0\| +  \int_{0}^{t} \|\varphi(s)\| d s ,
 \quad t > 0 .  
\end{equation} 

When using a two-level scheme with weights for an approximate solution of the problem (\ref{4.8}), (\ref{4.9}), (\ref{6.4}) 
instead of the equation (\ref{5.6}), we use the scheme
\begin{equation}\label{6.6}
 B \frac{y^{n+1} - y^{n}}{\tau } + \sum_{i=1}^{m} a_i A y_i^{n+\sigma} + C y^{n+\sigma} = \varphi^{n+\sigma } .
\end{equation} 
The stability estimate for the difference scheme  (\ref{5.7}), (\ref{5.8}), (\ref{6.6}) follows from the estimate (\ref{5.3}) and has the form
\[
 \Big (\|y^{n+1}\|_B^2 + \sum_{i=1}^{m} a_i \|y_i^{n+1}\|_{A}^2 \Big)^{1/2} \leq \|u_0\|_B + \sum_{k=0}^{n}\tau \|\varphi^{n+\sigma }\|_{B^{-1}} ,
 \quad n = 0,1, \ldots .
\] 
In the case of the scheme (\ref{5.7}), (\ref{6.6}), the equation (\ref{5.10}) is used to calculate $y_i^{n+1}, \ i = 1,2,\ldots, m,$ 
and for $y^{n+1}$ the equation 
\[
 \big (B + \sigma \tau (\mu A + C) \big ) y^{n+1} = \chi^{n} 
\] 
is solved with the same coefficient $\mu$. 
For the right-hand side, we have
\[
 \quad \chi^{n} = \big (B - (1-\sigma) \tau C \big ) y^{n} + \tau \varphi^{n+\sigma } - 
 \tau \sum_{i=1}^{m} a_i A \big ( (1-\sigma) y_i^{n} + \sigma \chi_i^{n} \big ) .
\] 
As before, the increase in the computational complexity of the approximate solution of the nonlocal problem (\ref{2.10}), (\ref{6.1}) is not of a fundamental nature.

\section{Numerical examples }\label{sec:7} 

The possibilities of the proposed computational algorithms will be illustrated by the results of the numerical solution of a model two-dimensional problem.
We will assume that the computational domain is unit square:
\[
 \Omega = \{ \bm x  \ | \ \bm x = (x_1,x_2), \ 0 < x_\alpha  < 1 , \ \alpha  = 1,2 \} ,
\]
with boundary $\partial \Omega$.
We define the operator 
\[
 \mathcal{A} w = - \triangle w, 
 \quad \bm x \in \Omega .
\] 

We consider the problem of relaxation of the initial state of the system under study when
\[
 f(\bm x,t) = 0,
 \quad u_0(\bm x) = x_1 x_2 \sin(\pi x_1) \sin(\pi x_2) .
\] 
In the numerical results presented below, the kernel is the stretching exponential function:
\[
 k(t) = \exp\big(- t^{\beta} \big),
 \quad 0 < \beta  < 1 . 
\] 

To numerically solve the problem (\ref{2.4}), (\ref{2.5}), we will use the standard
difference approximations in space \cite{SamarskiiTheory}.
We will introduce in the region $\Omega$ a uniform rectangular grid
\[
\overline{\omega}  = \{ \bm{x} \ | \ \bm{x} =\left(x_1, x_2\right), \quad x_\alpha  =
i_\alpha  h_\alpha , \quad i_\alpha  = 0,1,...,N_\alpha ,
\quad N_\alpha  h_\alpha  = 1 , \quad  \alpha  = 1,2 \} ,
\]
where $\overline{\omega} = \omega \cup \partial \omega$, 
$\omega$ is the set of internal mesh nodes, and $\partial \omega$ is the set of boundary mesh nodes.
For grid functions $w(\bm x)$ such that $w(\bm x) = 0, \ \bm x \notin \omega$, we define the Hilbert space
$H = L_2 (\omega)$, in which the dot product and norm are
\[
(w, u) = \sum_{\bm x \in  \omega} w(\bm{x}) u(\bm{x}) h_1 h_2,  \quad 
\| w \| =  (w,w)^{1/2}.
\]
For $w(\bm x) = 0, \ \bm x \notin \omega$, we define the grid Laplace operator $- A$ on the usual five-point pattern:
\[
  \begin{split}
  A w = & -
  \frac{1}{h_1^2} (w(x_1+h_1,x_2) - 2 w(\bm{x}) + w(x_1-h_1,x_2)) \\ 
  & - \frac{1}{h_2^2} (w(x_1,x_2+h_2) - 2 w(\bm{x}) + w(x_1,x_2-h_2)), 
  \quad \bm{x} \in \omega . 
 \end{split} 
\] 
For this grid operator (see, for example, \cite{SamarskiiTheory}), we have the property (\ref{2.8}).
On sufficiently smooth functions, the operator $A$ approximates the differential operator
$\mathcal{A}$ with an error
$\mathcal{O} \left(|h|^2\right)$, $|h|^2 = h_1^2+h_2^2$. 

\begin{figure}
\centering
\includegraphics[width=0.5\linewidth]{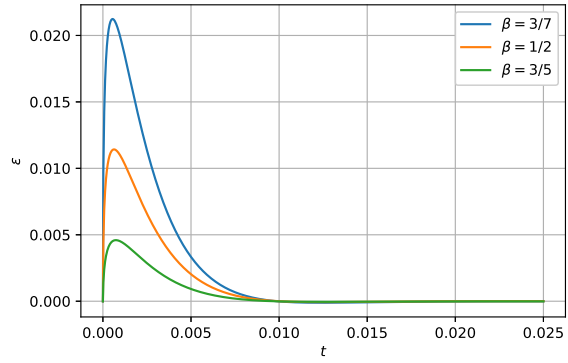}
\caption{Approximation error with $m = 12$ for stretched exponent.} 
\label{f-1}
\end{figure}

The construction of approximations (\ref{3.1}) is an independent task.
Considering that this problem is not the main one in our study, we will use the known results.
In the article \cite{mauro2018prony}, the stretched exponent approximations are obtained by the sum of exponentials for different $m$ under the additional constraint $\sum_{i=1}^{m}a_i = 1$.
In particular, the coefficients $a_i, b_i, \ i =1,2,\ldots, m,$ are given for $m = 8$ and $m = 12$.
We limited ourselves to the $m = 12$ option (see Table~1).
The approximation error $\varepsilon (t) = \accentset{\sim }{k}(t) - k(t)$ is most significant near $t = 0$ --- see  Fig.\ref{f-1}. 

\begin{center}
\begin{table}[htp]
\label{tab-1}
\caption{Parameters approximation with $m=12$ for stretched exponent (data from \cite{mauro2018prony})}
\centering
\begin{tabular}{|rrrrrrr|}
\hline
       &  &  $\beta = $ 3/7~~~  & & $\beta = $ 1/2~~~   &   & $\beta = $ 3/5~~~    \\
\hline
  $i$  & $a_i$~~~   & $b_i$~~~   &  $a_i$~~~    &  $b_i$~~~    &  $a_i$~~~    &  $b_i$~~~   \\
\hline
    1   &    0.02792  &  0.03816      &      0.01694  &  0.06265     &       0.01043  &  0.12022	\\
    2   &    0.09567  &  0.10117      &      0.08574  &  0.13381     &       0.08117  &  0.20610	\\
    3   &    0.13049  &  0.22822      &      0.14468  &  0.26816     &       0.17168  &  0.35680	\\
    4   &    0.13388  &  0.47142      &      0.15870  &  0.52050     &       0.19624  &  0.63293	\\
    5   &    0.12456  &  0.94243      &      0.14514  &  1.00410     &       0.16742  &  1.15481	\\
    6   &    0.10976  &  1.88828      &      0.12095  &  1.96395     &       0.12467  &  2.17404	\\
    7   &    0.09256  &  3.86312      &      0.09512  &  3.94401     &       0.08711  &  4.23811	\\
    8   &    0.07525  &  8.15604      &      0.07188  &  8.20241     &       0.05896  &  8.59467	\\
    9   &    0.05938  &  17.92388     &      0.05275  &  17.81155    &       0.03913  &  18.25401	\\
    10  &    0.04587  &  41.47225     &      0.03791  &  40.85894    &       0.02559  &  41.07522	\\
    11  &    0.03588  &  104.13591    &      0.02749  &  102.07104   &       0.01688  &  100.99297	\\
    12  &    0.06877  &  402.71691    &      0.04270  &  383.52267   &       0.02071  &  363.84147	\\
\hline
\end{tabular}
\end{table}
\end{center}

\begin{figure}
\centering
\includegraphics[width=0.5\linewidth]{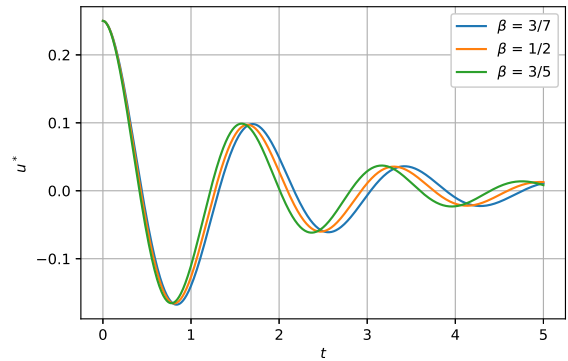}
\caption{Solution $u^* = u(\bm x^*, t)$ at the point $\bm x^* = (0.5,0.5)$ for different values of the parameter $\beta$.} 
\label{f-2}
\end{figure}

\begin{figure}
\centering
\begin{minipage}{0.49\linewidth}
\centering
\includegraphics[width=\linewidth]{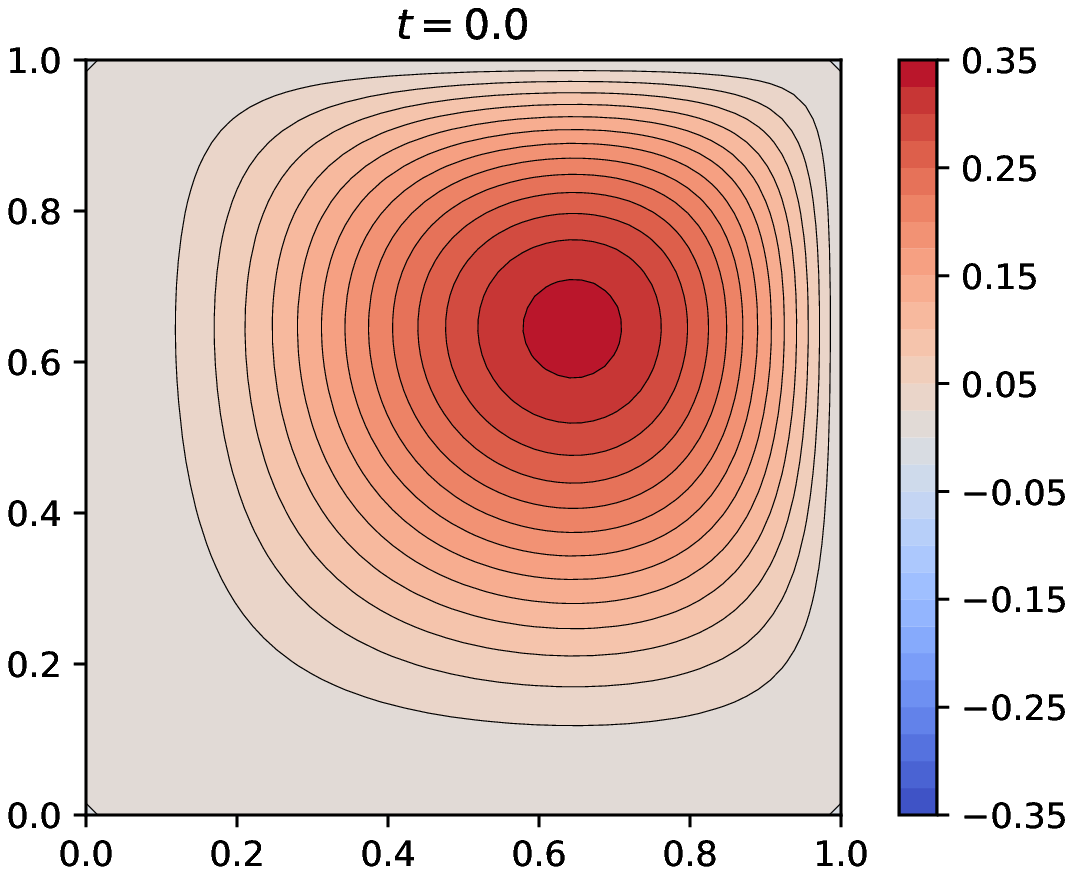} \\
\includegraphics[width=\linewidth]{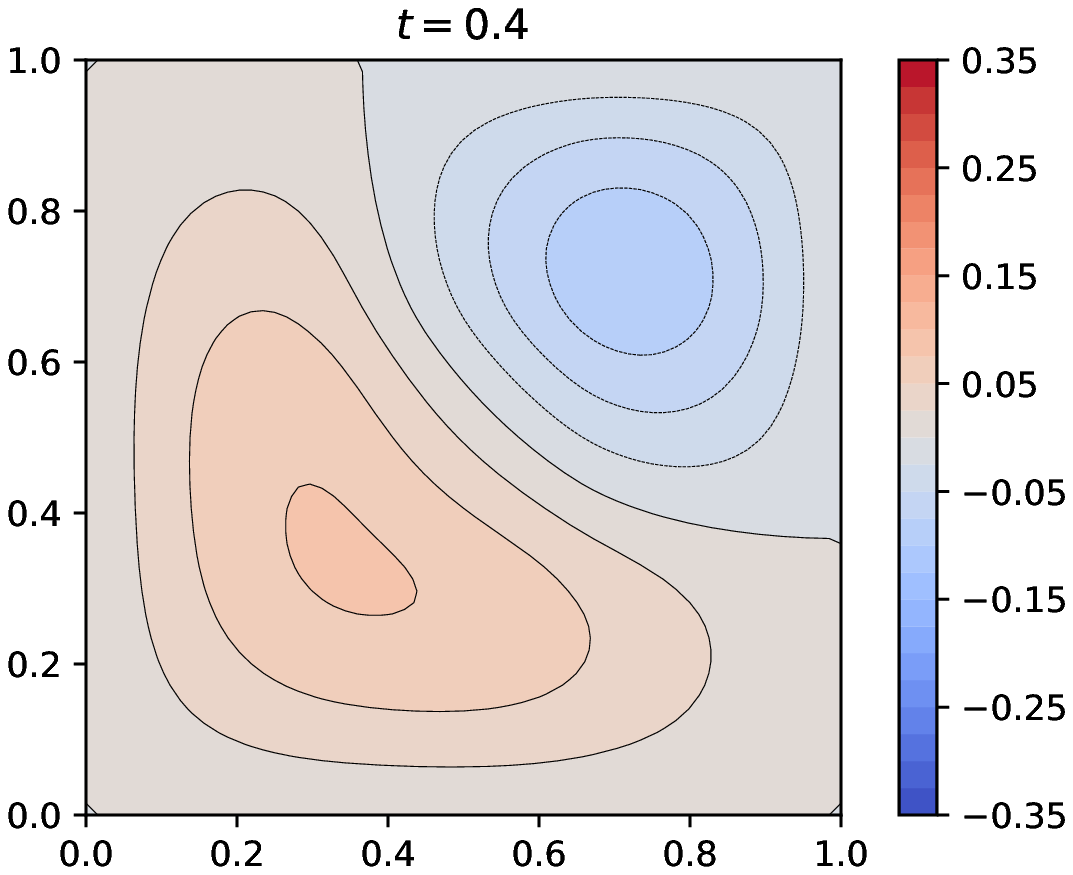} \\
\end{minipage}
\begin{minipage}{0.49\linewidth}
\centering
\includegraphics[width=\linewidth]{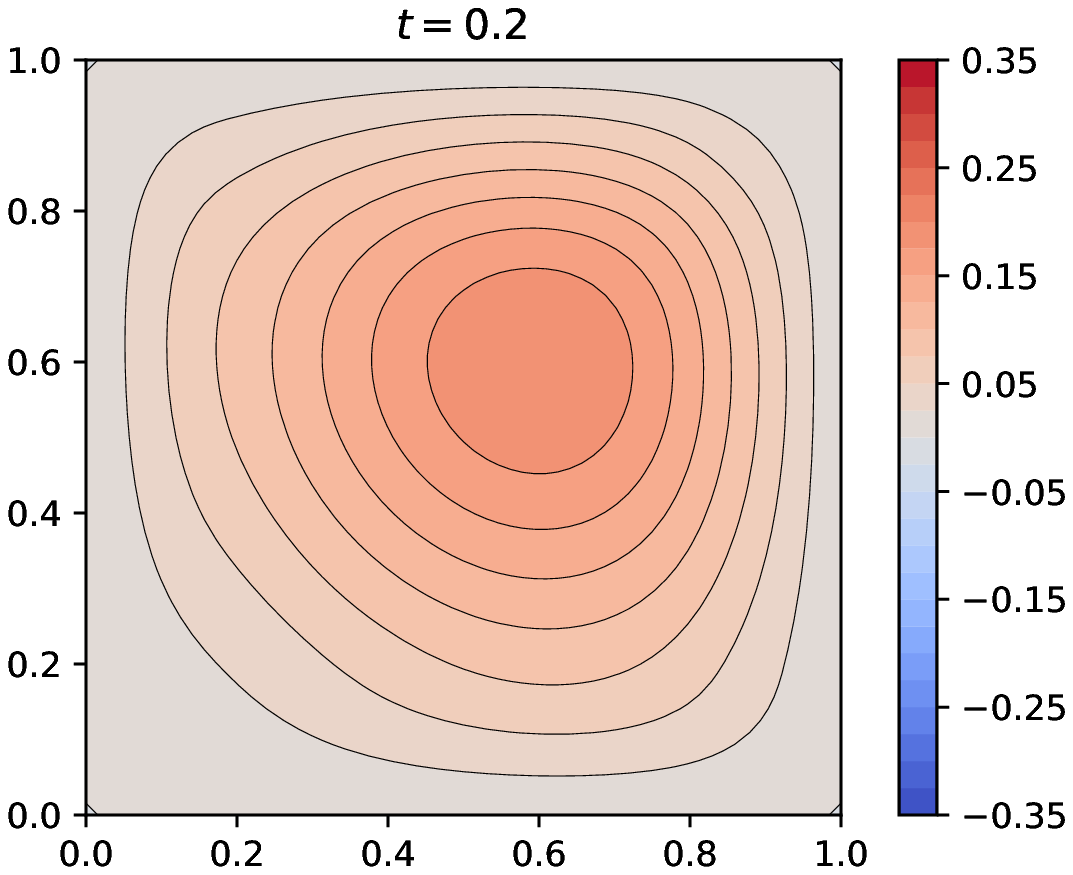} \\
\includegraphics[width=\linewidth]{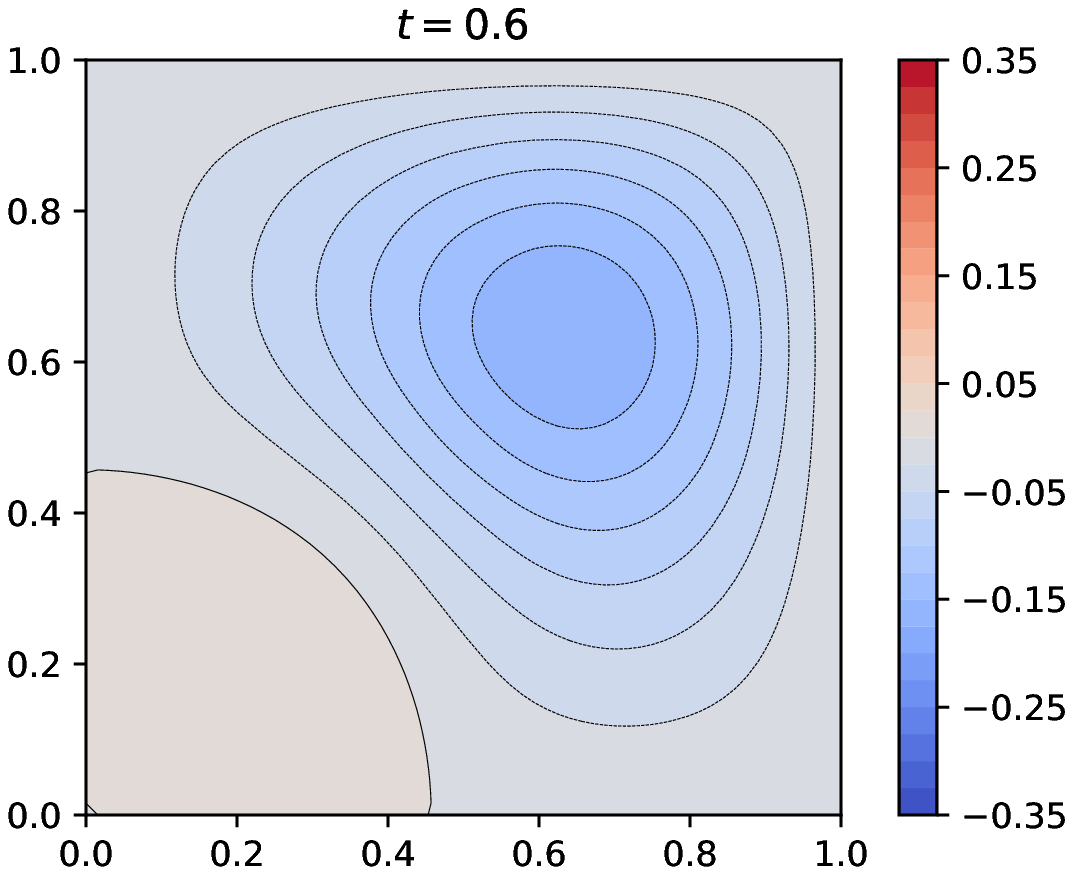} \\
\end{minipage}
\caption{Solution of the problem at separate points in time for $\beta = 1/2$.}
\label{f-3}
\end{figure}

The calculations were performed on the spatial grid $N_1 = N_2 = 64$. 
The approximate solution at different time steps was compared with the solution on a detailed grid with $N = 1000$, which was obtained using a second-order approximation scheme ($\sigma = 0.5$ in the scheme (\ref{5.6})--(\ref{5.8})). 
Fig.\ref{f-2} shows such a solution in the center of the computational domain.
We observe a pronounced oscillatory regime with amplitude damping.
We associate oscillations with hyperbolic properties, which are inherited by the considered equation with memory, and dissipation with parabolic properties. The solution at separate points in time for the problem with $\beta = 1/2$ is shown in  Fig.\ref{f-3}.

\begin{figure}
\centering
\includegraphics[width=0.45\linewidth]{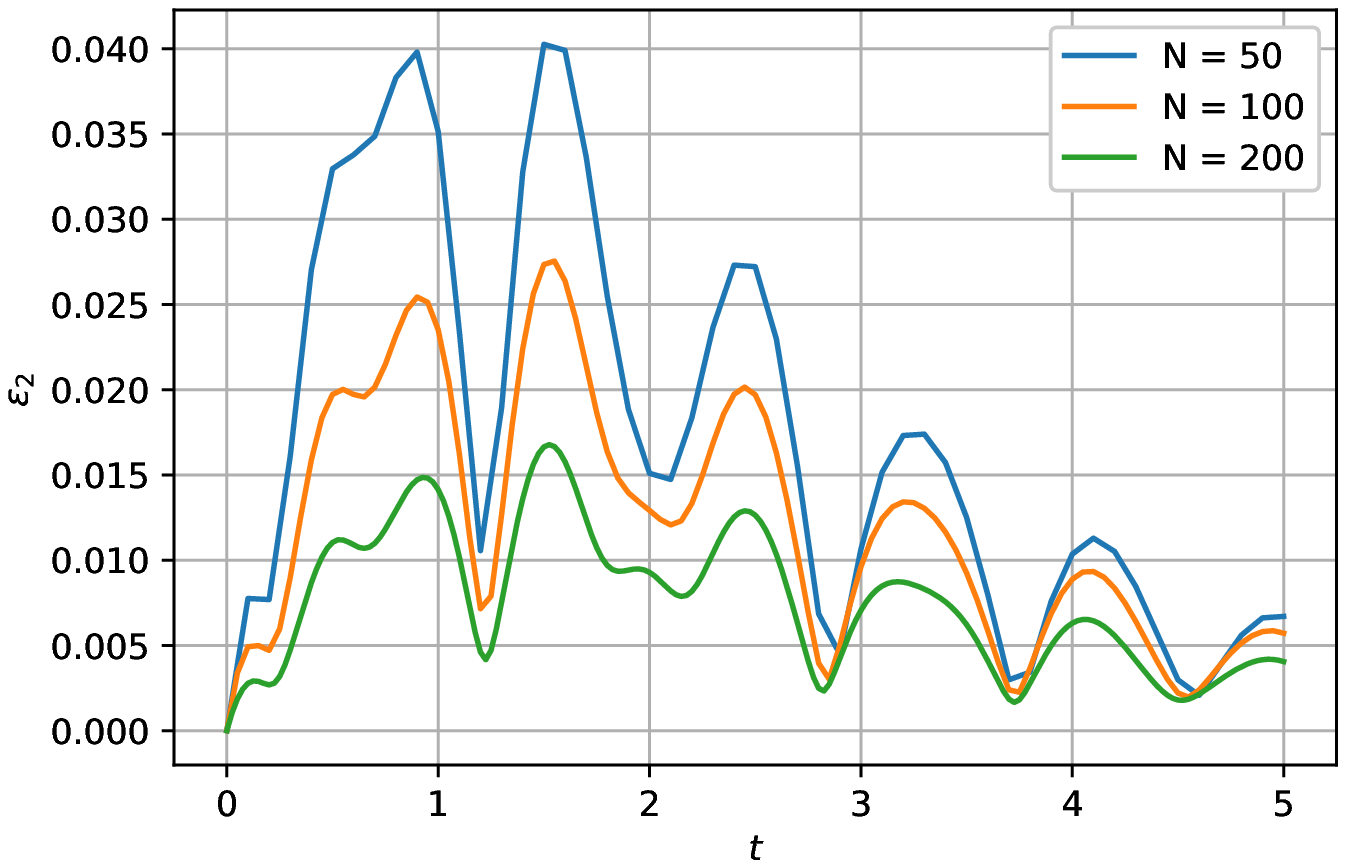} 
\includegraphics[width=0.45\linewidth]{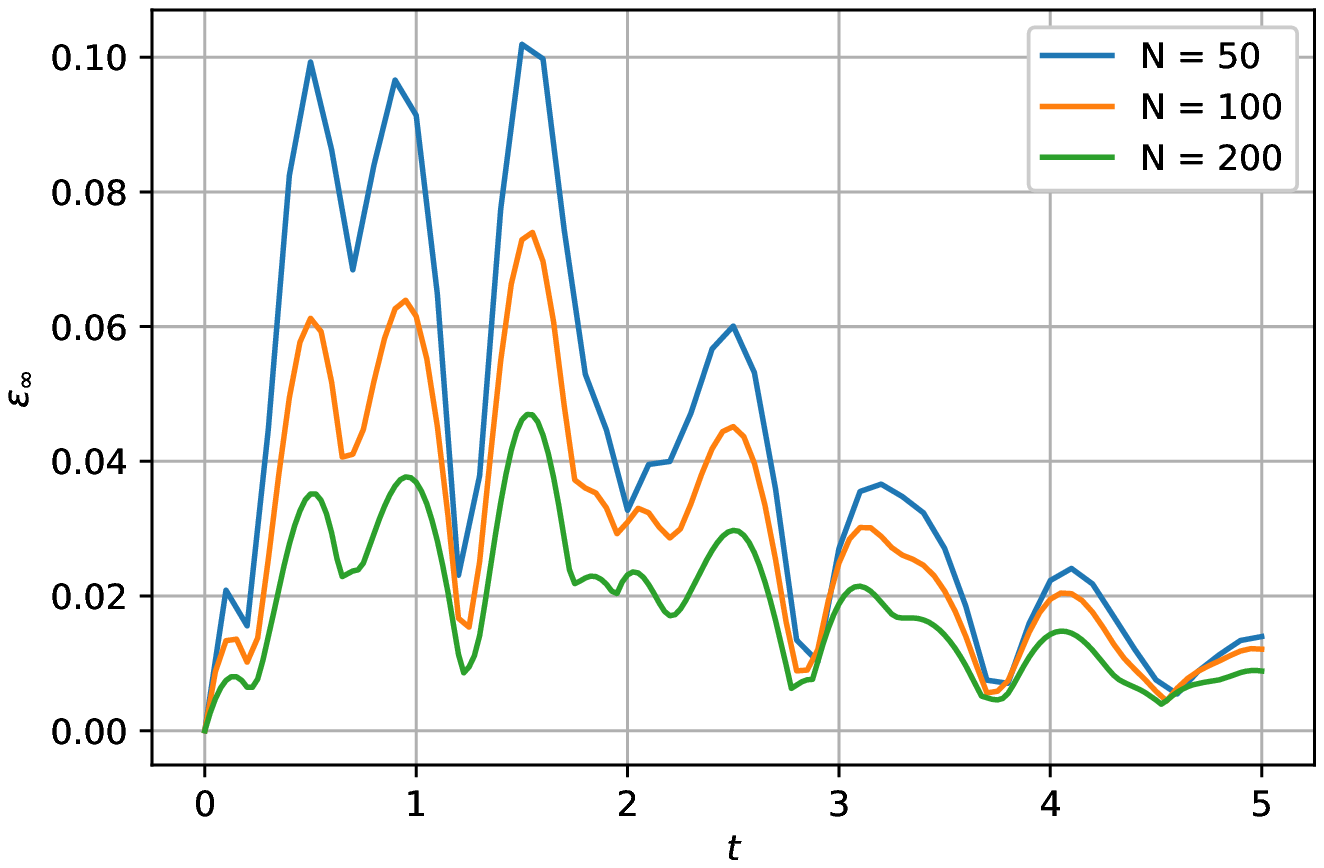}
\caption{Accuracy the implicit Euler scheme ($\sigma = 1$) for the problem with $\beta = 1/2$.}
\label{f-4}
\end{figure}

\begin{figure}
\centering
\includegraphics[width=0.45\linewidth]{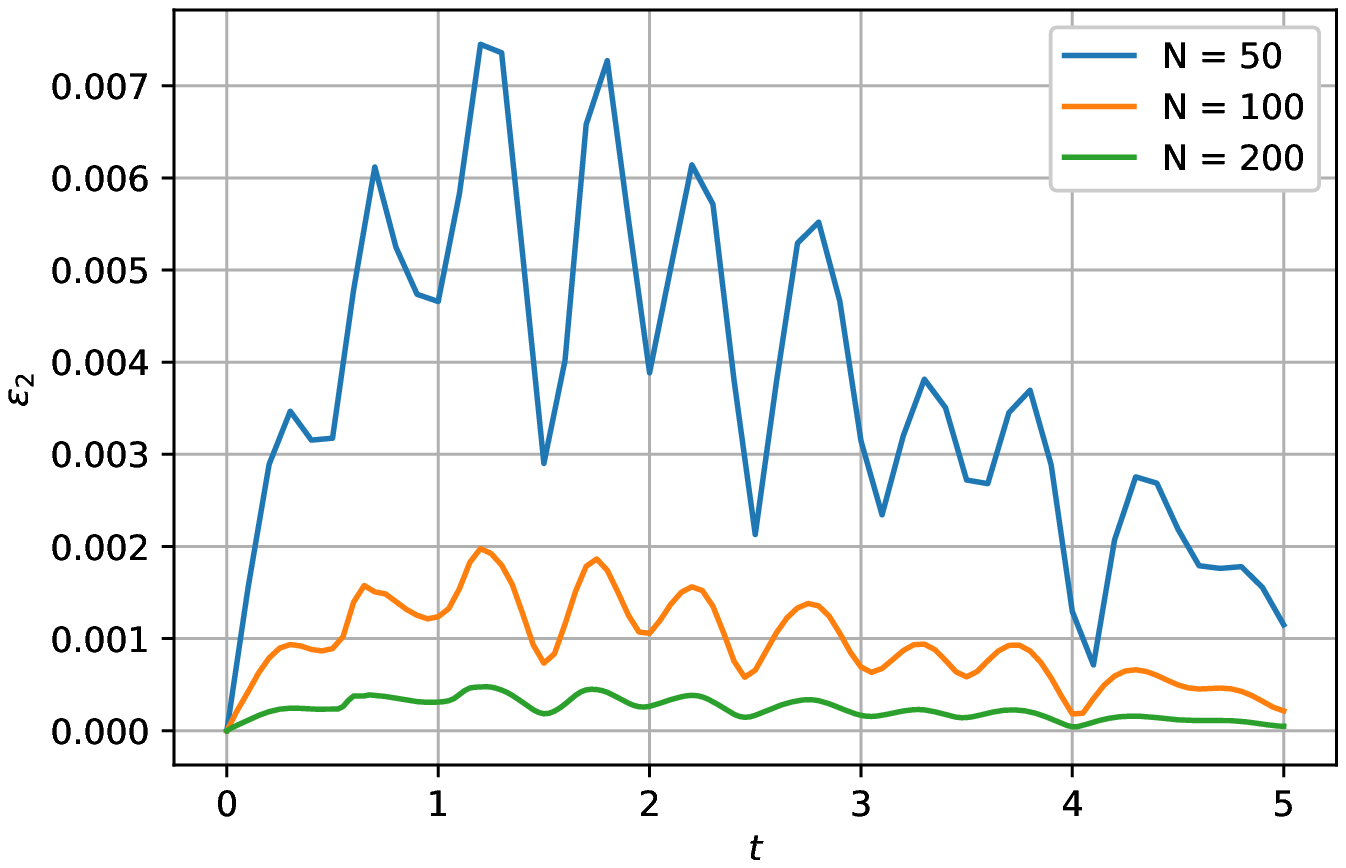} 
\includegraphics[width=0.45\linewidth]{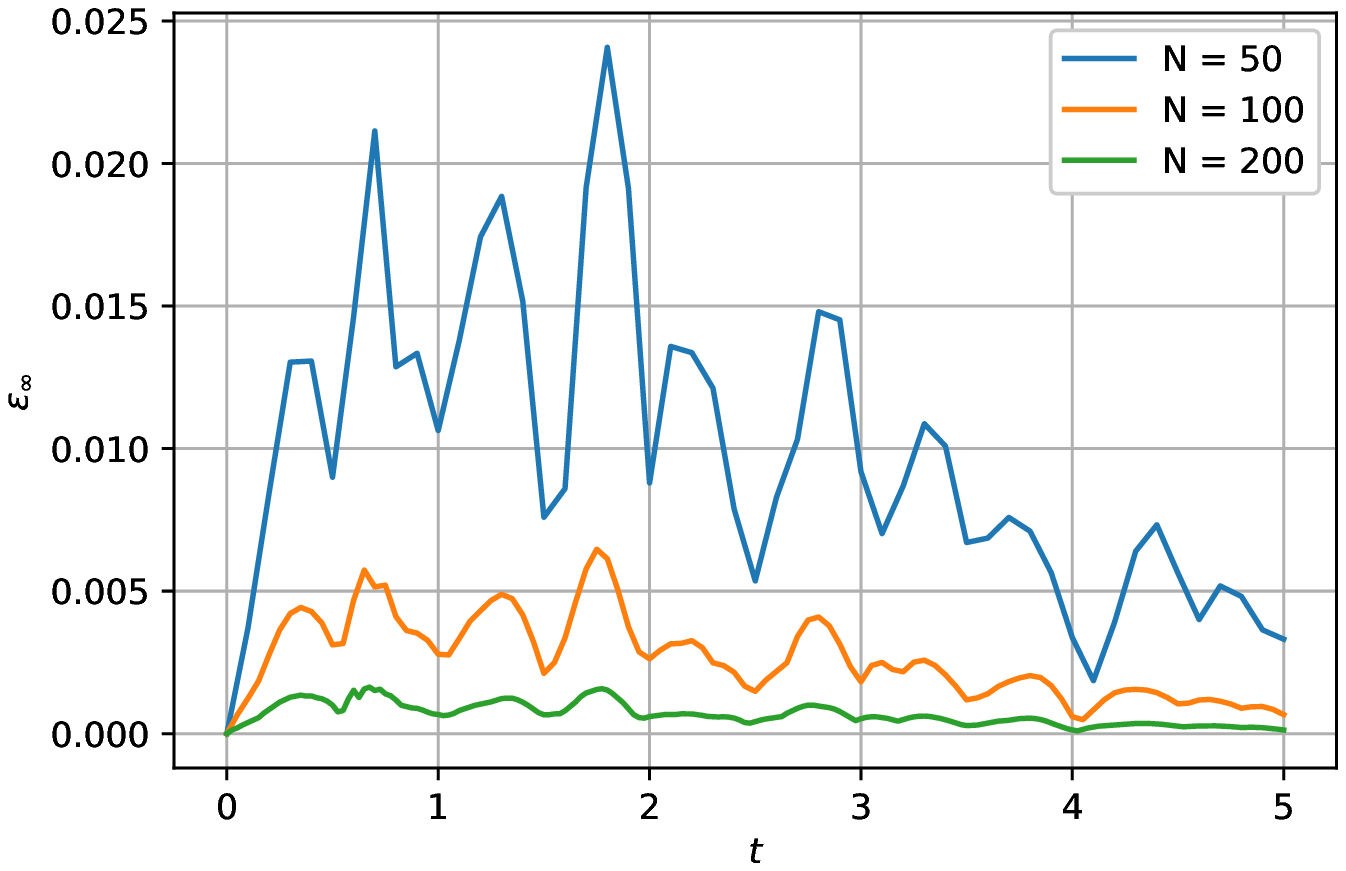} 
\caption{Accuracy the symmetric scheme ($\sigma = 0.5$) for the problem with $\beta = 1/2$.}
\label{f-5}
\end{figure}

The accuracy of the solution of problems with memory is estimated by the absolute discrepancy at individual points in time:
\[
 \varepsilon_2(t^n) = \|y (\bm x,t^n) - \bar{y}^n (\bm x,t^n)\| ,
 \quad  \varepsilon_\infty (t^n)  = \max_{\bm x \in \omega} |y (\bm x,t^n) - \bar{y}^n(\bm x)| ,
 \quad n = 0, \ldots, N ,
\] 
where $\bar{y}$ is the reference solution.
The accuracy when using the implicit Euler scheme is shown in Fig.\ref{f-4}.
We have given similar data for the symmetric scheme in Fig.\ref{f-5}. 
The calculated data are consistent with the above theoretical considerations for the accuracy of two-level scheme with weight (\ref{5.6})--(\ref{5.8}).

\section{Conclusions}\label{sec:8} 

\begin{enumerate}
 \item We considered the Cauchy problem for a first-order integrodifferential equation with a difference kernel and a self-adjoint positive definite operator in a finite-dimensional Hilbert space.
The stability of the solution concerning the right-hand side and initial conditions takes place under the assumption of a positive definite kernel.
The problems of the numerical solution of such problems are mainly associated with the need to operate with the solution for all previous moments of time.
 \item For an approximate solution to the stated evolutionary problem with memory, a well-known approach is used, based on the approximation of the difference kernel by the sum of exponentials --- the Prony series.
In this case, we pass from a nonlocal in time problem to a local problem for a system of weakly coupled evolution equations.
The nonlocality of the processes under study is associated with additional ordinary differential equations for auxiliary functions.
A priori estimates for the solution of the Cauchy problem for this system of evolution equations are proved.
 \item In the numerical solution, we use standard two-level time approximations.
The unconditional stability of two-level schemes with weights is proved under standard weight constraints.
The transition to a new level in time is provided by solving the usual problem for an approximate solution and explicitly recalculating auxiliary functions.
 \item Our results are generalized to more general evolutionary problems with memory.
We supplemented the theoretical consideration with data from the numerical solution of a model two-dimensional problem, when the kernel is the stretching exponential function.
\end{enumerate}


\end{document}